\documentclass[12pt]{article}
\title{Minimality of hyperplane arrangements and basis of 
local system cohomology}
\author{Masahiko Yoshinaga\thanks{
Department of Mathematics, Kyoto University, 
Kyoto 606-8502, Japan, 
email: mhyo@math.kyoto-u.ac.jp
}}
\date{\today}

\usepackage{amsmath, amssymb, amsthm, verbatim} 
\usepackage{eucal} 

\newtheorem{Def}{Definition}[section]
\newtheorem{Prop}[Def]{Proposition}
\newtheorem{Thm}[Def]{Theorem}
\newtheorem{Lemma}[Def]{Lemma}
\newtheorem{Conj}[Def]{Conjecture}
\newtheorem{Cor}[Def]{Corollary}
\newtheorem{Rem}[Def]{Remark}
\newtheorem{Problem}[Def]{Problem}
\newtheorem{Example}[Def]{Example}

\newcommand{\bbC}{\mathbb{C}}

\newcommand{\bbP}{\mathbb{P}}

\newcommand{\bbR}{\mathbb{R}}

\newcommand{\bbZ}{\mathbb{Z}}

\newcommand{\calA}{\mathcal{A}}

\newcommand{\calF}{\mathcal{F}}

\newcommand{\calL}{\mathcal{L}}

\newcommand{\A}{\calA}
\newcommand{\bch}{\mathsf{bch}}

\newcommand{\ch}{\mathsf{ch}}

\newcommand{\F}{\calF}

\newcommand{\sfD}{\mathsf{D}}

\newcommand{\M}{\mathsf{M}}

\newcommand{\uch}{\mathsf{uch}}

\makeatletter
\def\codim{\mathop{\operator@font codim}\nolimits}
\def\Coker{\mathop{\operator@font Coker}\nolimits}
\def\grad{\mathop{\operator@font grad}\nolimits}
\def\Hom{\mathop{\operator@font Hom}\nolimits}
\def\Im{\mathop{\operator@font Im}\nolimits}
\def\Ker{\mathop{\operator@font Ker}\nolimits}
\def\Loc{\mathop{\operator@font Loc}\nolimits}
\def\rank{\mathop{\operator@font rank}\nolimits}
\def\Rep{\mathop{\operator@font Rep}\nolimits}
\def\Sep{\mathop{\operator@font Sep}\nolimits}
\def\tr{\mathop{\operator@font tr}\nolimits}
\def\sign{\mathop{\operator@font sign}\nolimits}
\makeatother

\begin{document}
\maketitle

\begin{abstract}
The purpose of this paper is applying 
minimality of hyperplane arrangements to 
local system cohomology groups. 
It is well known that twisted cohomology groups 
with coefficients in a generic rank one local system 
vanish except in the top degree, and bounded chambers 
form a basis of the remaining cohomology group. 
We determine precisely when this phenomenon happens 
for two-dimensional arrangements. 
\end{abstract}

\section{Introduction}
\label{sec:intro}

The purpose of this paper is applying 
minimality of hyperplane arrangements to 
local system cohomology groups. In 
\S\ref{subsec:min} and \S\ref{subsec:nonres}, 
we will recall basic notions and results on these topics. 
In \S\ref{subsec:plan}, we will give the plan of 
the paper.

\subsection{Minimality of hyperplane arrangements}
\label{subsec:min}

Let $\calA=\{H_1, \dots, H_n\}$ be a hyperplane 
arrangement in $\bbC^\ell$. Namely a finite set 
of affine hyperplanes. We assume each hyperplane 
$H_i=\{\alpha_i=0\}\subset\bbC^\ell$ is defined by 
an affine linear equation $\alpha_i$. 
We denote the complement of hyperplanes 
by $\M(\A)=\bbC^\ell\setminus\bigcup_{i=1}^n H_i$. 

After the discovery of combinatorial description 
of the cohomology ring $H^*(\M(\A), \bbZ)$ 
\cite{orl-sol} and $K(\pi, 1)$-property 
for simplicial arrangements \cite{del-simp}, 
it has been revealed that the complement $\M(\A)$ of 
a hyperplane arrangement $\A$ has a very special 
homotopy type 
among other complex affine varieties. 
Especially, the following minimality 
seems one of the most peculiar properties 
to $\M(\A)$ \cite{dim-pap, ran-mor, ps-h, fal-hom}.

\begin{Thm}
\normalfont\label{thm:min}
(Minimality of arrangements.) 
The complement $\M(\A)$ is homotopy equivalent to 
a finite minimal CW-complex $X$. Namely, $X$ satisfies 
the following minimality: The number of $k$-dimensional 
cells $\sharp\{k\mbox{-}\dim\mbox{ cells}\}$ is equal to 
the $k$-th Betti number $b_k(X)$. 
\end{Thm}

The minimality is expected to 
be useful for computations of local system 
cohomology groups. 
An immediate corollary is 
the following upper bounds for dimensions 
of rank one local system cohomology groups, 
which were conjectured by Aomoto and first 
proved in \cite{coh-mor} by using another method. 
\begin{Cor}
\normalfont\label{cor:bound}
Let $\calL$ be a complex rank one local 
system on $\M(\A)$. Then the dimension of 
$\calL$-coefficients cohomology group is 
bounded by Betti number:
$$
\dim H^k(\M(\A), \calL)\leq b_k(\M(\A)), 
$$
for $k=0, 1, \dots, \ell$. 
\end{Cor}

For further applications of the minimality 
to computations of local system cohomology groups, 
the description 
of the minimal CW-complex $X$, in particular 
the attaching map of each cell, is needed. 
However Theorem \ref{thm:min} does not tell it. 
It should be noted that the proof of Theorem \ref{thm:min} 
is based on Morse theoretic arguments. The constructions of 
cells are relying on a transcendental method, 
namely using gradient flows of a Morse function. 

Both of the two recent approaches to the 
problem of describing attaching maps of minimal cells 
are: 
\begin{itemize}
\item assuming $\A$ is defined over the real numbers $\bbR$, and 
\item describing attaching maps by using combinatorial 
structure of chambers. 
\end{itemize}
However they used different methods. 
\begin{itemize}
\item In \cite{yos-lef}, we studied 
Lefschetz's hyperplane section theorem for 
$\M(\A)$, and described the attaching maps of 
the top cells. 
\item In \cite{sal-sett}, Salvetti and Settepanella 
developed discrete Morse theory on the Salvetti 
complex, and then described the minimal cell complex by 
using discrete Morse flows. 
\end{itemize}
See \cite{del-min, del-sett} for subsequent 
developments. Furthermore, in \cite{gai-sal}, $2$-dimensional 
algebraic minimal chain complex is described. 
The present article can be considered as 
a counterpart of \cite{gai-sal}.

\subsection{Non-resonant local systems}
\label{subsec:nonres}

A nonempty intersection of elements of $\A$ is 
called an {\em edge}. We denote by $L(\A)$ 
the set of edges. An edge $X\in L(\A)$ is called 
a {\em dense edge} if the localization 
$\A_X=\{H\in\A\mid H\supset X\}$ is indecomposable. 
We denote by $\sfD(\A)\subset L(\A)$ the set of 
dense edges. 

Let $\lambda=(\lambda_1, \dots, \lambda_n)\in\bbC^n$. 
Then $\lambda$ determines a rank one representation of 
$\pi_1(\M(\A))$ by 
$\pi_1(\M(\A))\ni\gamma\longmapsto\exp
(\int_\gamma\sum_{i=1}^n\lambda_id\log\alpha_i)
\in\bbC^*$ and the associated local system 
$\calL=\calL_\lambda$. In other words, 
$\calL$ is determined by the local monodromy 
$q_i=e^{2\pi\sqrt{-1}\lambda_i}\in\bbC^*$ around 
each hyperplane $H_i$. For an edge $X\in L(\A)$, 
denote $q_X=\prod_{X\subset H_i}q_i$. 
We also denote the half 
twist by $q_i^{1/2}=e^{\pi\sqrt{-1}\lambda_i}$. 

We can embed the affine space $\bbC^\ell$ in 
$\bbC\bbP^\ell$ as 
$\bbC^\ell=\bbC\bbP^\ell\setminus H_\infty$. 
We call $\A_\infty:=\{\overline{H}\mid H\in\A\}\cup
\{H_\infty\}$ the projective closure of $\A$. 
The monodromy of $\calL_\lambda$ around the hyperplane 
at infinity $H_\infty$ is $\prod_{i=1}^n q_i^{-1}$. 
It is natural to define 
$q_\infty=\prod_{i=1}^n q_i^{-1}$. 

The structure of the cohomology group 
$H^k(\M(\A), \calL)$ with local system coefficients 
has been studied well 
\cite{aom, esv, koh, stv}. In particular, it is known that 
if $\calL$ is generic, then the cohomology 
vanishes except in $k=\ell$. Among others, let us recall two 
results in this direction. 
(\cite{dt-det, lib-eig, cdo})

\begin{Thm}
\normalfont\label{thm:dt}
(\cite{dt-det})
Suppose that $\A$ is defined over $\bbR$ and 
the local system $\calL_\lambda$ 
satisfies 
\begin{equation}
\label{eq:dt}
q_X\neq 1\mbox{, for }\forall X\in\sfD(\A_\infty). 
\end{equation}
Then 
\begin{equation}
\label{eq:gen}
H^k(\M(\A), \calL_\lambda)=
\left\{
\begin{array}{cl}
0 & \mbox{for } k\neq \ell, \\
&\\
\bigoplus\limits_{C\in\bch(\A)}\bbC\cdot [C],& 
\mbox{ for } k=\ell, 
\end{array}
\right.
\end{equation}
where $\bch(\A)$ stands for the set of all 
bounded chambers. 
A chamber $[C]$ can be considered as 
a locally finite cycle, in other words, 
an element of Borel-Moore homology 
$[C]\in H_\ell^{BM}(\M(\A))$. 
In (\ref{eq:gen}) we identify the chamber 
$C$ with cohomology via the canonical isomorphism 
$H_\ell^{BM}(\M(\A))\simeq H^\ell(\M(\A))$. 
\end{Thm}

\begin{Def}
\normalfont\label{def:infty}
$\sfD_\infty(\A_\infty):=\{X\in\sfD(\A_\infty)\mid 
X\subset H_\infty\}$. 
\end{Def}

\begin{Thm}
\normalfont\label{thm:cdo}
(\cite{cdo})
Suppose that the local system $\calL_\lambda$ 
satisfies 
\begin{equation}
\label{eq:cdo}
q_X\neq 1\mbox{, for }\forall X\in\sfD_\infty(\A_\infty). 
\end{equation}
Then 
$$
H^k(\M(\A), \calL_\lambda)\simeq
\left\{
\begin{array}{cl}
0 & \mbox{for } k\neq \ell, \\
&\\
\bbC^{|\chi(\M(\A))|}& 
\mbox{ for } k=\ell, 
\end{array}
\right.
$$
where $\chi(\M(\A))$ is the Euler characteristic 
of $\M(\A)$. 
\end{Thm}

\subsection{Plan of the paper}
\label{subsec:plan}

The purpose of this paper is to refine 
vanishing results 
Theorem \ref{thm:dt} and Theorem \ref{thm:cdo} for $\ell=2$ 
by using minimal complex arising from 
minimal CW-decomposition of $\M(\A)$. 
We will prove that 
the assertion (\ref{eq:gen}) of 
Theorem \ref{thm:dt} is true under the weaker 
assumption (\ref{eq:cdo}). Furthermore, 
if $\A$ is indecomposable, we also prove 
that the assumption can not be weakened 
any more. 
Our main result asserts that 
(\ref{eq:cdo}) and (\ref{eq:gen}) are 
equivalent. (For $\ell=2$.) 

In \S\ref{sec:flag}, we treat combinatorial 
structures of chambers, which will play a crucial role 
in the study of minimal complex. 

In \S\ref{sec:min}, we will describe the minimal 
cochain complex arising from Lefschetz's hyperplane 
section theorem. Particularly, we treat the 
case $\ell=2$ in details. 

In \S\ref{sec:appl}, we prove the main result, that is, 
for an indecomposable two dimensional arrangement 
$\A$, conditions (\ref{eq:cdo}) and (\ref{eq:gen}) 
are equivalent. 

\section{Chambers and flags}
\label{sec:flag}


\subsection{Involution on unbounded chambers}
\label{subsec:invol}

Let $\A$ be a hyperplane arrangement in $\bbR^\ell$. 
We denote the set of chambers, bounded chambers, 
unbounded chambers by $\ch(\A), \bch(\A), \uch(\A)$, 
respectively. Note that $\ch(\A)=\bch(\A)\sqcup\uch(\A)$. 

Let $C\in\uch(\A)$ be an unbounded chamber. 
Then the closure $cl(C)$ in the projective space 
$\bbR\bbP^\ell$ intersects the hyperplane 
$H_\infty$ at infinity. 
\begin{Def}
\normalfont 
Let $C\in\uch(\A)$. 
(i) Define $X(C)$ to be the smallest 
subspace of $H_\infty$ which contains 
$cl(C)\cap H_\infty$. 
(ii) There exists a unique chamber which is the 
opposite with respect to $cl(C)\cap H_\infty$. 
We denote the opposite chamber by $C^\lor$ (see 
Figure \ref{fig:invol}). 
Obviously we have $C^{\lor\lor}=C$. 
\end{Def}
See Figure \ref{fig:invol} for an example. 
In this figure, $X(C_1)=X(C_4)=H_\infty$, and 
$X(C_2)=X(C_3)=cl(C_2)\cap H_\infty$.

\begin{figure}[htbp]
\begin{picture}(100,200)(0,0)
\thicklines

\put(200,40){\circle*{5}}
\put(200,40){\line(-1,-1){30}}
\put(200,40){\line(0,-1){40}}
\put(200,40){\line(1,-1){30}}
\put(200,40){\line(0,1){160}}

\qbezier(200,40)(230,70)(230,80)
\put(230,80){\line(0,1){120}}
\qbezier(200,40)(170,70)(170,80)
\put(170,80){\line(0,1){120}}

\put(100,40){\line(1,0){180}}
\qbezier(280,40)(300,40)(300,60)
\put(300,60){\line(0,1){140}}

\put(100,110){\line(1,0){207}}
\put(100,90){\line(5,1){207}}

\put(80,35){$H_\infty$}

\put(140,80){$C_1$}
\put(180,80){$C_2$}
\put(210,80){$C_3$}
\put(250,80){$C_4$}

\put(140,140){$C_4^\lor$}
\put(180,140){$C_2^\lor$}
\put(210,140){$C_3^\lor$}
\put(250,140){$C_1^\lor$}

\put(170,30){\vector(3,1){25}}
\put(110,25){$cl(C_2)\cap H_\infty$}
\put(302,110){\line(0,-1){50}}
\qbezier(302,60)(302,38)(280,38)
\put(280,38){\line(-1,0){76}}
\put(330,90){\vector(-3,-1){25}}
\put(333,87){$cl(C_4)\cap H_\infty$}

\put(208,5){$C_2^\lor$}
\put(178,5){$C_3^\lor$}
\put(310,30){$C_4^\lor$}

\end{picture}
\caption{$C$ and $C^\lor$}
\label{fig:invol}
\end{figure}

\begin{Def}
\normalfont
Define the involution $\iota$ by 
\begin{eqnarray*}
\iota:\uch(\A)&\longrightarrow&\uch(\A)\\
C&\longmapsto&C^\lor
\end{eqnarray*}
\end{Def}

We now characterize dense edges contained in 
$H_\infty$ by using $X(C)$. First we prove an 
easy lemma. 

\begin{Lemma}
\normalfont
Let $\A$ be an essential central arrangement 
in $\bbR^\ell$. Then the following are equivalent. 
\begin{itemize}
\item[(1)] $\A$ is indecomposable. 
\item[(2)] There exist $H\in\A$ and 
$C\in\ch(\A)$ such that $cl(C)\cap H=\{0\}$. 
\item[(3)] For any $H\in\A$, there exists 
$C\in\ch(\A)$ 
such that $cl(C)\cap H=\{0\}$. 
\end{itemize}
\end{Lemma}

\proof
Let $H\in\A$ and consider the deconing 
${\bf d}_H\A$ with respect to $H$. 
Note that ${\bf d}_H\A$ is an 
affine arrangement of rank $(\ell-1)$. 
Using \cite[\S3.3]{ot-hyp}, 
$\A$ is indecomposable if and only if 
the $\beta$-invariant of ${\bf d}_H\A$ is nonzero. 
By the famous result of Zaslavsky 
\cite{zas-face}, it is equivalent to 
the existence of bounded chambers of the deconing 
${\bf d}_H\A$. Choose a bounded chamber 
of ${\bf d}_H\A$, and 
let $C$ be its cone. 
Then $cl(C)\cap H=\{0\}$. This proves 
(1) $\Rightarrow$ (3). The other implications can also be 
similarly proved. 
\qed

Using the above lemma, we obtain the following. 

\begin{Prop}
\label{prop:dense}
\normalfont
Let $\A$ be an affine arrangement in $\bbR^\ell$. 
An edge $X\in L(\A_\infty)$ satisfies 
$X\in\sfD_\infty(\A_\infty)$ if and only if 
$X=X(C)$ for some $C\in\uch(\A)$. 
\end{Prop}

\subsection{Generic flags}
\label{subsec:flag}

Let $\F$ be a generic flag 
 in $\Bbb{R} ^\ell$ 
$$
\F : 
\emptyset = \F ^{-1} \subset 
\F ^0 
\subset 
\F ^1 
\subset 
\cdots 
\subset 
\F ^\ell = \Bbb{R} ^\ell ,  
$$  
where  
each $\F ^q$ 
is a {\it generic} $q$-dimensional affine subspace, 
that is, 
$\dim \F ^q \cap X = q+ \dim X-\ell$ 
for $X \in L (\A_\infty)$. 
Let 
$\{ h_1 , \ldots h_\ell \}$ 
be a system of defining equations of $\F$, 
that is, 
$$
\F ^{q} 
= \{ h_{q+1} = \cdots = h_\ell =0 \} ,  \mbox{ for }
q=0, 1, \dots, \ell-1, 
$$  
where each $h_i$ is an affine linear form on $\Bbb{R} ^\ell$. 
Using the flag $\F$, we decompose the set of chambers 
into several subsets. 
\begin{Def}
\normalfont
Define 
$$ 
\ch ^q ( \A ) 
= \{ 
C \in \ch (\A) \mid 
C \cap \F ^q \not= \emptyset \quad \mbox{and} \quad 
C \cap \F ^{q-1} = \emptyset 
\}, 
$$ 
for $q=0, 1, \dots, \ell$. 
\end{Def}
\begin{Prop}
\normalfont(\cite{yos-lef}) 
\label{prop:betti}
$\sharp\ch^q(\A)=b_q(\M(\A))$. 
\qed
\end{Prop}

\begin{Rem}
\normalfont
\label{rem:zas}
The above proposition gives a refinement 
of Zaslavsky's formula $\sum_{i=0}^\ell b_i(\M(\A))=
\sharp\ch(\A)$, (\cite{zas-face}). 
\end{Rem}

We assume that $\F$ satisfies the following : 
For $q=0, \ldots ,\ell$, 
$\F ^{q} _{>0}$ 
 denotes  
$$ 
\{ h_{q+1}= h_{q+2} = \cdots = h_\ell =0, h_{q} >0 
\} . 
$$  
\begin{enumerate} 
\item For an arbitrary chamber $C$, 
if belonging to  
$\ch ^{q} (\A)$, then 
$C \cap \F ^{q} \subset \F ^{q} _{>0}$.     
\item For any two 
$X$, $X^\prime \in L(\A)$ with 
$\dim X=\dim X^\prime = \ell- q$ (i.e. 
satisfying 
$X \cap \F^{q} = \{ pt \}$ 
and 
$X ^\prime \cap \F^{q} = \{ pt \}$),  
if $X \not= X^\prime$, 
$$ 
h_{q} (X \cap \F^{q}) 
\not= 
h_{q} (X ^\prime \cap \F^{q})  . 
$$ 
\end{enumerate} 
In the remainder of the paper we fix a 
generic flag $\F$ satisfying the above 
conditions. And also fix the orientation 
of $\F^q$ by the oriented basis 
$(\partial_{h_1}, \dots, \partial_{h_q})$ 
of the tangent space $T_x\F^q$. 

Next we further decompose $\ch^q(\A)$ into two subsets. 
\begin{Def}
\normalfont
Define subsets $\bch^q(\A)$ and $\uch^q(\A)$ of 
$\ch^q(\A)$ by 
\begin{eqnarray*}
\bch^q(\A)&=&\{C\in\ch^q(\A)\mid C\cap\F^q
\mbox{ is bounded}\}, \\
\uch^q(\A)&=&\{C\in\ch^q(\A)\mid C\cap\F^q
\mbox{ is unbounded}\}. 
\end{eqnarray*}
We note that $\bch^\ell(\A)=\bch(\A)$. 
\end{Def}

\begin{Example}
\normalfont
Let us consider the arrangement of 
four lines $\A=\{H_1, H_2, H_3, H_4\}$ with 
a generic flag $\F$ as in 
Figure \ref{fig:buch}. 

\begin{figure}[htbp]
\begin{picture}(100,130)(0,0)
\thicklines

\put(273,3){$H_4$}
\put(233,3){$H_3$}
\put(105,3){$H_2$}
\put(65,3){$H_1$}

\multiput(120,10)(-40,0){2}{\line(1,1){120}}
\multiput(230,10)(40,0){2}{\line(-1,1){120}}

\put(50,100){$\F^2=\Bbb{R}^2$}

\put(60,35){$\F^0$}
\put(60,30){\circle*{4}}
\put(325,35){$\F^1$}
\multiput(50,30)(5,0){60}{\circle*{2}}

\put(215,85){$C_0^\lor$}
\put(192.5,62.5){$C_3$}
\put(170,40){$C_2$}

\put(150,60){$C_1$}
\put(172.5,82.5){$C_4$}
\put(195,105){$C_1^\lor$}

\put(130,80){$C_0$}
\put(152.5,102.5){$C_3^\lor$}
\put(175,125){$C_2^\lor$}



\end{picture}
\caption{$\bch^q(\A)$ and $\uch^q(\A)$.}
\label{fig:buch}
\end{figure}
Then we have by definition 
$$
\begin{array}{lclcl}
\ch^0(\A)=\{C_0\}, &&
\ch^1(\A)=\{C_1, C_2, C_3, C_0^\lor\}, &&
\ch^2(\A)=\{C_1^\lor, C_2^\lor, C_3^\lor, C_4\},\\
\bch^0(\A)=\{C_0\}, &&
\bch^1(\A)=\{C_1, C_2, C_3\}, &&
\bch^2(\A)=\{C_4\}, \\
\uch^0(\A)=\emptyset, &&
\uch^1(\A)=\{C_0^\lor\}, &&
\uch^2(\A)=\{C_1^\lor, C_2^\lor, C_3^\lor\}.
\end{array}
$$
\end{Example}
\begin{Thm}
\normalfont\label{thm:invol}
The involution $\iota$ induces 
a bijection 
$$
\iota:\bch^{q-1}(\A)\stackrel{\sim}{\longrightarrow}
\uch^q(\A). 
$$
\end{Thm}

\proof
Suppose $C\in\bch^{q-1}(\A)$, that is, 
$C\cap\F^{q-1}$ is bounded. 
Then $C^\lor\cap\F^{q-1}=\emptyset$. 
By the assumption on 
the flag, $C\cap\F^q$ is unbounded. Since $\F^q$ is 
generic, 
$cl(\F^q)$ intersects $cl(C)\cap H_\infty$ transversally. 
Hence $C^\lor\cap\F^q\neq\emptyset$ and unbounded. 
We have $C^\lor\in\uch^q(\A)$. Conversely if 
$C^\lor\in\uch^q(\A)$, then $C^{\lor\lor}=C$ intersects 
$\F^{q-1}$. Suppose $C\cap\F^{q-1}$ is unbounded. In this 
case, $C^\lor$ also intersects $\F^{q-1}$. This 
contradicts the fact 
$C^\lor\in\uch^q(\A)\subset\ch^q(\A)$. 
\qed

\begin{Cor}
\normalfont
\label{cor:bij}
$\sharp\bch^{q-1}(\A)=\sharp\uch^q(\A)$. \qed
\end{Cor}

\begin{Rem}
\normalfont
(1) Corollary \ref{cor:bij} together with 
Proposition \ref{prop:betti} and 
$\sharp\ch^q(\A)=\sharp\bch^q(\A) + \sharp\uch^q(\A)$, gives a 
``bijective proof'' for Zaslavsky's formula 
$\sharp\bch(\A)=\sum_{i=0}^\ell (-1)^{\ell-i}b_i(\M(\A))$. 


(2) The bijective correspondence 
(Theorem \ref{thm:invol}) plays a crucial role 
in \S\ref{sec:appl}. 
\end{Rem}

\section{Minimal complexes}
\label{sec:min}

Let $\A$ be an essential real arrangement and 
$\F$ be a generic flag as in the previous section. 
Set $F=\F^{\ell-1}\otimes\bbC$ the complexification of 
$\F^{\ell-1}$. 
Compare the complexified complement 
$\M(\A)$ with the generic hyperplane 
section $\M(\A)\cap F$. 
Lefschetz's hyperplane section theorem \cite{le-ham} 
tells us that $\M(\A)$ is homotopy equivalent to 
the space obtained from $\M(\A)\cap F$ by 
attaching some $\ell$-dimensional cells. 
Namely we have the following homotopy 
equivalence: 
$$
\M(\A)\approx(\M(\A)\cap F)\cup_{\varphi_i}
\bigcup_i D^\ell, 
$$
where $\varphi_i:\partial D^\ell\longrightarrow 
\M(\A)\cap F$ is the attaching map. 
In \cite{yos-lef}, we described the 
homotopy type of the attaching maps. 
The $\ell$-dimensional cells are naturally encoded by 
the set $\ch^\ell(\A)$ of chambers which 
do not intersect $\F^{\ell-1}$. 
By using the description of attaching maps, 
we constructed a cochain complex 
$$
(\bbC[\ch^q(\A)], d_\calL)_{q=0}^\ell : 
\cdots\longrightarrow
\bbC[\ch^q(\A)]
\stackrel{d_\calL}{\longrightarrow}
\bbC[\ch^{q+1}(\A)]
\longrightarrow
\cdots
$$
which computes local system cohomology 
groups for arbitrary rank one local system $\calL$. Namely, we have 
$
H^*
(\bbC[\ch^\bullet(\A)], d_\calL)
\simeq
H^*
(\M(\A), \calL)$. 
In \S\ref{subsec:lefmin}, we shall describe 
the cochain complex 
$(\bbC[\ch^\bullet(\A)], d_\calL)$ 
based on \cite{yos-lef}, and in 
\S\ref{subsec:2dim} 
we investigate the case $\ell=2$ closely. 

\subsection{Minimal complex arising from Lefschetz's Theorem}
\label{subsec:lefmin}

\begin{Def}
\normalfont
(Separating hyperplanes) Let $C_1, C_2\in\ch(\A)$ be 
chambers. Define 
$$
\Sep(C_1, C_2)=\{H\in\A\mid
H\mbox{ separates $C_1$ and $C_2$}\}. 
$$
And also 
$$
q_{\Sep(C_1, C_2)}^{1/2}=
\prod_{H_i\in\Sep(C_1, C_2)}q_i^{1/2}. 
$$
\end{Def}
To describe the coboundary map 
$d_\calL:\bbC[\ch^q(\A)]\rightarrow\bbC[\ch^{q+1}(\A)]$, 
we need the notion of degree map 
$$
\deg:\ch^q(\A)\times
\ch^{q+1}(\A)\longrightarrow\bbZ, 
$$
which we will define below. 

Suppose $C\in\ch^q(\A)$ and $C'\in\ch^{q+1}(\A)$ 
are given. 
Let $D=D^q\subset\F^q$ be a $q$-dimensional 
ball with sufficiently large radius so that 
every $0$-dimensional edge $x\in L(\A\cap\F^q)$ 
is in the interior of $D^q$. 
There exists a tangent vector field 
$U(x)\in T_x\F^q$ for 
$x\in D$ 
which satisfies the following properties: 
\begin{itemize}
\item if $x\in\partial D$, then 
$U(x)\notin T_x(\partial D)$, and $U(x)$ 
directs inside of $D$, 
\item if $x\in H$ with $H\in\A$, then 
$U(x)\notin T_x(H\cap\F^q)\subset T_x\F^q$ and 
$U(x)$ directs the side in which $C'$ is contained. 
\end{itemize}
From the properties, we have $U(x)\neq 0$ for $x\in 
\partial(cl(C)\cap D)$, where $cl(C)$ is the 
closure of $C$ in $\F^q$. Roughly speaking, the degree 
$\deg(C, C')$ is defined to be the degree of the Gauss 
map
$$
\frac{U}{|U|}: 
\partial(cl(C)\cap D)\longrightarrow S^{q-1}. 
$$
\begin{Def}
\normalfont
Let $C\in\ch^q(\A)$ and $C'\in\ch^{q+1}(\A)$. 
Fix $U$ as above. Then define $\deg(C, C')$ as 
follows. 
\begin{itemize}
\item[(0)] When $q=0$, then $\deg (C, C')=1$. 
\item[(1)] When $q=1$, then $cl(C)\cap D\simeq [-1, 1]$. 
In this case $S^0\simeq\{\pm 1\}$. 
The degree of the Gauss maps 
$g:=\frac{U}{|U|}:\{\pm 1\}\longrightarrow \{\pm 1\}$ 
is defined by 
$$
\deg(g)=\left\{
\begin{array}{cl}
0  &   \mbox{if $g(\{\pm 1\})=\{+1\}$ or $g(\{\pm 1\})=\{-1\}$}, \\
1  &   \mbox{if $g(\pm 1)=\pm 1$}, \\
-1 &   \mbox{if $g(\pm 1)=\mp 1$}. 
\end{array}
\right.
$$
\item[(2)] When $q\geq 2$, 
$$
\deg(C, C')=\deg
\left(
\frac{U}{|U|}: 
\partial(cl(C)\cap D)\longrightarrow S^{q-1}
\right). 
$$
(It is easily seen that $\deg(C, C')$ does not depend on $U$.) 
\end{itemize}
\end{Def}
Now let us define the map 
$$
d_\calL:
\bbC[\ch^q(\A)]\longrightarrow
\bbC[\ch^{q+1}(\A)]
$$
by
\begin{equation}
\label{eq:d}
\ch^q(\A)\ni [C]\longmapsto
\sum_{C'\in\ch^{q+1}(\A)}
\deg(C, C')\cdot
\left(q_{\Sep(C, C')}^{1/2}-
q_{\Sep(C, C')}^{-1/2}\right)
\cdot[C']. 
\end{equation}

\begin{Thm}
\normalfont
(\cite[6.4.1]{yos-lef})
With notation as above, 
$(\bbC[\ch^\bullet(\A)], d_\calL)$ is a 
cochain complex. Furthermore, 
$$
H^*(\bbC[\ch^\bullet(\A)], d_\calL)\simeq
H^*(\M(\A), \calL). 
$$
\end{Thm}


In the above formula (\ref{eq:d}), 
the degree $\deg(C, C')\in\bbZ$ is difficult to determine. 
The author wonders how to compute 
$\deg(C, C')$. Let us pose a problem which 
might be interesting from the view point of 
combinatorics of polytopes. 

\begin{Problem}
\normalfont
Let $P\subset\bbR^d$ be a bounded $d$-dimensional 
convex polytope. Let 
$\{F_e\}_{e\in E}$ be the set of 
facets (i.e., $(d-1)$-dimensional faces). 
Let $U(x)\in T_x\bbR^d$ be a vector field on $\bbR^d$. 
Suppose that $U$ satisfies $U(x)\neq 0$ when 
$x\in\partial P$ and, furthermore, $U(x)\notin T_x F_e$ 
for any point $x\in F_e$ in a facet. 
We can associate a sign vector 
$X\in\{+1, -1\}^E$ by 
$$
X(e)=\left\{
\begin{array}{ll}
+1 & \mbox{if $U$ directs outside of $P$ on $F_e$}, \\
-1 & \mbox{if $U$ directs inside of $P$ on $F_e$}. 
\end{array}
\right.
$$
Then how to compute the degree 
$\deg\left(\frac{U}{|U|}:\partial P
\rightarrow S^{d-1}\right)$ of the Gauss map 
from the sign vector $X\in\{\pm1\}^E$? 
\end{Problem}

\subsection{The case $\ell=2$}
\label{subsec:2dim}

In this section, we look at the minimal complex 
$(\bbC[\ch^\bullet(\A)], d_\calL)$ for 
$\ell=2$ more closely. 

First note that $\ch^0(\A)=\{C_0\}$ consists of a 
chamber. The map 
$d_\calL: \bbC[\ch^0(\A)]\longrightarrow\bbC[\ch^1(\A)]$ 
is determined by $d_\calL([C_0])$, which is 
$$
d([C_0])= 
\sum_{C\in\ch^1(\A)}
\left(
q_{\Sep(C_0, C)}^{1/2}-
q_{\Sep(C_0, C)}^{-1/2}
\right)\cdot [C]. 
$$
As in \S\ref{subsec:invol}, we decompose 
$\ch^1(\A)=\bch^1(\A)\sqcup\uch^1(\A)$. Note that 
by Theorem \ref{thm:invol}, 
$\uch^1(\A)=\{C_0^\lor\}$ consists of a chamber 
which is the opposite one of $C_0$. 
The second coboundary map 
$d_\calL:\bbC[\ch^1(\A)]\longrightarrow
\bbC[\ch^2(\A)]$ is given by the formula 
(\ref{eq:d}). The degree 
$\deg(C, C')$ behaves differently according as 
$C\in\bch^1(\A)$ or $C\in\uch^1(\A)$. 
\begin{itemize}
\item[(i)] Suppose $C\in\bch^1(\A)$. Then $C\cap\F^1$ 
is a closed interval, the boundary (two points) 
can be expressed as 
$(H\cap\F^1)\cup(H'\cap\F^1)$ for 
$H, H'\in\A$. 
$\deg(C, C')$ can be computed as 
$$
\deg(C, C')=
\left\{
\begin{array}{cl}
1& \mbox{ if } H, H'\in\Sep(C, C'), \\
-1& \mbox{ if } H, H'\notin\Sep(C, C'), \\
0& \mbox{ others.}
\end{array}
\right.
$$
\item[(ii)] Suppose $C\in\uch^1(\A)$. Then $C\cap\F^1$ 
is an unbounded interval, the boundary (a point) 
can be expressed as 
$H\cap\F^1$. 
$\deg(C, C')$ can be computed as 
$$
\deg(C, C')=
\left\{
\begin{array}{cl}
-1& \mbox{ if } H\notin\Sep(C, C'), \\
0& \mbox{ if } H\in\Sep(C, C'). 
\end{array}
\right.
$$
\end{itemize}
In particular, we have, 
\begin{Lemma}
\label{lem:deg}
\normalfont
Let $C\in\bch^1(\A)$. The boundary of $C\cap\F^1$ is 
expressed as $(H\cap\F^1)\cup(H'\cap\F^1)$. Then 
$$
\deg(C, C^\lor)=
\left\{
\begin{array}{cl}
1& \mbox{ if $H$ and $H'$ are not parallel,}\\
-1& \mbox{ if $H$ and $H'$ are parallel.}
\end{array}
\right.
$$
\end{Lemma}

\begin{Example}
\normalfont
\label{ex:cpx}
Consider the arrangement of 
four lines $\A=\{H_1, H_2, H_3, H_4\}$ in $\bbR^2$ and 
a generic flag $\F$ as in Figure \ref{fig:cpx}. 
\begin{figure}[htbp]
\begin{picture}(100,130)(0,0)
\thicklines

\put(95,0){$H_1$}
\put(175,0){$H_2$}
\put(215,0){$H_3$}
\put(335,0){$H_4$}

\put(100,10){\line(2,1){210}}
\put(340,10){\line(-2,1){210}}
\multiput(180,10)(40,0){2}{\line(0,1){120}}

\put(50,100){$\F^2=\Bbb{R}^2$}

\multiput(50,15)(5,0){66}{\circle*{2}}
\put(355,19){$\F^1$}

\put(60,15){\circle*{4}}
\put(60,19){$\F^0$}

\put(140,65){$C_0$}
\put(160,25){$C_1$}
\put(198,25){$C_2$}
\put(240,25){$C_3$}
\put(260,65){$C_0^\lor$}
\put(160,110){$C_3^\lor$}
\put(198,110){$C_2^\lor$}
\put(240,110){$C_1^\lor$}

\put(190,65){$D$}

\end{picture}
\caption{Example \ref{ex:cpx}.}
\label{fig:cpx}
\end{figure}
Then 
$$
\begin{array}{lclcl}
\bch^0(\A)=\{C_0\}& & 
\bch^1(\A)=\{C_1, C_2, C_3\} & & 
\bch^2(\A)=\{D\}\\
\uch^0(\A)=\emptyset & & 
\uch^1(\A)=\{C_0^\lor\}& & 
\uch^2(\A)=\{C_1^\lor, C_2^\lor, C_3^\lor\}. 
\end{array}
$$
The coboundary map $d_\calL:\bbC[\ch^0]\rightarrow\bbC[\ch^1]$ 
is determined by 
$$
d_\calL([C_0])=
(q_1^{\frac{1}{2}}-q_1^{-\frac{1}{2}})[C_1]+
(q_{12}^{\frac{1}{2}}-q_{12}^{-\frac{1}{2}})[C_2]+
(q_{123}^{\frac{1}{2}}-q_{123}^{-\frac{1}{2}})[C_3]+
(q_{1234}^{\frac{1}{2}}-q_{1234}^{-\frac{1}{2}})[C_0^\lor], 
$$
and $d_\calL:\bbC[\ch^1]\rightarrow\bbC[\ch^2]$ is 
as follows. 
{\small 
$$
\begin{array}{llr}
d_\calL([C_1])=
\underline{(q_{1234}^{\frac{1}{2}}-q_{1234}^{-\frac{1}{2}})}[C_1^\lor]&
+(q_{124}^{\frac{1}{2}}-q_{124}^{-\frac{1}{2}})[C_2^\lor]
&
+(q_{12}^{\frac{1}{2}}-q_{12}^{-\frac{1}{2}})[D]
\\

d_\calL([C_2])=
&
-\underline{(q_{14}^{\frac{1}{2}}-q_{14}^{-\frac{1}{2}})}[C_2^\lor]
&
-(q_{1}^{\frac{1}{2}}-q_{1}^{-\frac{1}{2}})[D]
\\
d_\calL([C_3])= &
+(q_{134}^{\frac{1}{2}}-q_{134}^{-\frac{1}{2}})[C_2^\lor]
+\underline{(q_{1234}^{\frac{1}{2}}-q_{1234}^{-\frac{1}{2}})}[C_3^\lor]&
\\
d_\calL([C_0^\lor])=
-(q_{1}^{\frac{1}{2}}-q_{1}^{-\frac{1}{2}})[C_1^\lor]&
-(q_{13}^{\frac{1}{2}}-q_{13}^{-\frac{1}{2}})[C_2^\lor]
-(q_{123}^{\frac{1}{2}}-q_{123}^{-\frac{1}{2}})[C_3^\lor]&
\end{array}
$$}
\end{Example}

The coefficients of the diagonals have another 
expressions. 
Observe that $X(C_1)=X(C_3)=H_\infty$ and 
$X(C_2)=\overline{H_2}\cap\overline{H_3}\cap H_\infty$. 
Since $q_\infty=q_{1234}^{-1}$ and 
$q_{X(C_2)}=q_2q_3q_\infty=q_{14}^{-1}$, we have 
\begin{eqnarray*}
\underline{(q_{1234}^{\frac{1}{2}}-q_{1234}^{-\frac{1}{2}}
)}
&=&
-(q_{X(C_1)}^{\frac{1}{2}}-q_{X(C_1)}^{-\frac{1}{2}})
=-(q_{X(C_3)}^{\frac{1}{2}}-q_{X(C_3)}^{-\frac{1}{2}}), 
\\
-\underline{(q_{14}^{\frac{1}{2}}-q_{14}^{-\frac{1}{2}})}
&=&
q_{X(C_2)}^{\frac{1}{2}}-q_{X(C_2)}^{-\frac{1}{2}}. 
\end{eqnarray*}
In general, we have 
\begin{Prop}
\label{prop:coeff}
\normalfont
Let $C\in\bch^1(\A)$. 
Then the coefficient 
of $[C^\lor]$ in $d_\calL([C])$ is given by 
$\pm(q_{X(C)}^{1/2}-q_{X(C)}^{-1/2})$. 
\end{Prop}

\proof 
Let $H\in\A$. Then $H$ separates 
$C$ and $C^\lor$ if and only if 
$\overline{H}$ does go through 
$X(C)\in H_\infty$. 
Using $q_1q_2, \dots q_n q_\infty=1$, we have 
$$
q_{\Sep(C, C^\lor)}=q_{X(C)}^{-1}. 
$$
Hence 
$\pm(q_{\Sep(C, C^\lor)}^{1/2}-q_{\Sep(C, C^\lor)}^{-1/2})=
\mp(q_{X(C)}^{1/2}-q_{X(C)}^{-1/2})$. 
\qed

For use in the next section, 
we analyze the induced map 
$$
\bbC[\bch^1(\A)]
\hookrightarrow
\bbC[\ch^1(\A)]
\stackrel{d_\calL}{\longrightarrow}
\bbC[\ch^2(\A)]
\twoheadrightarrow
\bbC[\uch^2(\A)]. 
$$
Consider the composed map 
$\overline{d_\calL}: 
\bbC[\bch^1(\A)]
\longrightarrow
\bbC[\uch^2(\A)]$. 
As Theorem \ref{thm:invol}, the bases of 
the source and the target of $\overline{d_\calL}$ 
are naturally identified by the involution $\iota$. 
Thus the determinant $\det(\overline{d_\calL})\in
\bbC$ makes sense. 
The matrix $\overline{d_\calL}$ is expressed 
by an upper triangular matrix, and the 
determinant can be computed.

\begin{Thm}
\label{thm:det}
\normalfont
The determinant 
$\det(\overline{d_\calL})$ can be 
expressed as 
\begin{equation}
\label{eq:det}
\det(\overline{d_\calL})=\pm
\prod_{X\in \sfD_\infty(\A_\infty)}
(q_{X}^{1/2}-q_X^{-1/2})^{n_X}, 
\end{equation}
where $n_X$ is a positive integer. 
\end{Thm}

\proof
First note that, for $C\in\uch(\A)$, 
$X(C)$ is either $0$-dimensional or equal to 
$H_\infty$. 
We call an unbounded chamber $C\in\uch(\A)$ 
{\em narrow} (resp. {\em wide}) if 
$X(C)\subset H_\infty$ is $0$-dimensional 
(resp. $X(C)=H_\infty$). 
We decompose $\bbC[\bch^1(\A)]$ and 
$\bbC[\uch^2(\A)]$ into direct  sum of 
subspaces. 
Set 
\begin{eqnarray*}
N^1=\bbC[\{C\in\bch^1(\A)\mid C:\mbox{ narrow}\}], &
W^1=\bbC[\{C\in\bch^1(\A)\mid C:\mbox{ wide}\}], \\
N^2=\bbC[\{C\in\uch^2(\A)\mid C:\mbox{ narrow}\}], &
W^2=\bbC[\{C\in\uch^2(\A)\mid C:\mbox{ wide}\}]. \\
\end{eqnarray*}
Then clearly $\bbC[\bch^1(\A)]=W^1\oplus N^1$ and 
$\bbC[\uch^2(\A)]=W^2\oplus N^2$. The map 
$\overline{d_\calL}$ preserves $N^i$. 
Furthermore, the matrix presentation of 
$\overline{d_\calL}|_{N^1}:N^1\rightarrow N^2$ is 
diagonal. Indeed suppose that 
$C\in\bch^1(\A)$ is a narrow chamber 
with walls $H\cap\F^1$ and 
$H'\cap\F^1$. Then $H$ and $H'$ are 
parallel. By definition of degree map, 
$d_\calL([C])$ is a linear combination of 
chambers which are put between $H$ and $H'$. 
The opposite chamber $C^\lor$ is the unique such 
element in $\uch^2(\A)$. 
By Proposition \ref{prop:coeff}, we obtain 
the explicit formula 
$$
\overline{d_\calL}([C])=(q_{X(C)}^{1/2}-q_{X(C)}^{-1/2})
[C^\lor]
$$
for a narrow chamber $C\in\bch^1(\A)$. 
Next we consider $W^1$ and $W^2$. 
Since $\bbC[\bch^1]/N^1\simeq W^1$ and 
$\bbC[\uch^2]/N^2\simeq W^2$, we have the 
induced map $\widetilde{d_\calL}:W^1\rightarrow W^2$. 
This map is again expressed by a diagonal matrix. 
Indeed, for a wide chamber $C\in\bch^1(\A)$, we have 
\begin{eqnarray*}
\widetilde{d_\calL}([C])&=&
-(q_\infty^{1/2}-q_\infty^{-1/2})[C^\lor]\\
&=&
-(q_{X(C)}^{1/2}-q_{X(C)}^{-1/2})[C^\lor]. 
\end{eqnarray*}
Thus 
\begin{eqnarray*}
\det\left(
\overline{d_\calL}:\bbC[\bch^1]\rightarrow\bbC[\uch^2]
\right)&=&
\det(\overline{d_\calL}|_{N^1})\cdot
\det(\widetilde{d_\calL})\\
&=&
\pm\prod_{C\in\bch^1(\A)}(q_{X(C)}^{1/2}-q_{X(C)}^{-1/2}). 
\end{eqnarray*}
By Proposition \ref{prop:dense}, $X(C)$ in the above 
formulas runs all dense edges contained in 
$H_\infty$. Hence we obtain (\ref{eq:det}). 
\qed

\begin{Cor}
\normalfont
\label{cor:nonzero}
The map 
$\overline{d_\calL}: 
\bbC[\bch^1(\A)]
\longrightarrow
\bbC[\uch^2(\A)]$ 
is nondegenerate if and only if 
$q_X\neq 1$ for any dense edge 
$X\in\sfD_\infty(\A_\infty)$ in $H_\infty$. 
\end{Cor}

The decomposability of $\A$ is related to 
$W^2$ as follows. We omit the proof 
(cf. Figure \ref{fig:buch} and Figure \ref{fig:cpx}). 
\begin{Prop}
\normalfont
\label{prop:indec}
For $\ell=2$, 
$\A$ is decomposable if and only if $\dim W^2=1$. 
\end{Prop}

\section{An application}
\label{sec:appl}

As we saw in the previous sections, the basis 
of our cochain complex is encoded by the set 
of chambers. There is also an involution $\iota$ 
among unbounded chambers. In this section, we 
prove that if the monodromies around dense edges 
at infinity are nontrivial, then the bases 
corresponding to unbounded chambers $C$ and 
$C^\lor=\iota(C)$ are cancelled 
each other, and finally, only 
bounded chambers survive. This leads to 
a proof of the refined version of vanishing theorem. 

Our main result is the following. 
\begin{Thm}
\label{thm:main}
\normalfont
Let $\A$ be an indecomposable line 
arrangement in $\bbR^2$. Let $\calL$ be 
a rank one local system. Then the following 
are equivalent. 
\begin{itemize}
\item[(i)] $q_X\neq 1$ for any dense edge 
$X\in\sfD_\infty(\A_\infty)$ contained in $H_\infty$. 
\item[(ii)] 
$$
H^k(\M(\A), \calL)=\left\{
\begin{array}{cl}
0&\mbox{ for }k=0, 1, \\
&\\
\bigoplus\limits_{C\in\bch(\A)}\bbC\cdot[C]& 
\mbox{ for }k=2. 
\end{array}
\right.
$$
\item[(iii)] $H^2(\A, \calL)$ is generated 
by $\{[C]\mid C\in\bch(\A)\}$. 
\end{itemize}
\end{Thm}
\begin{Rem}
\normalfont
(i)$\Rightarrow$(ii)$\Rightarrow$(iii) holds for 
any arrangement $\A$ (without indecomposability). 
However (iii)$\Rightarrow$(i) requires the 
indecomposability of $\A$. (See Remark \ref{rem:counter}.) 
For comments to the higher dimensional cases ($\ell\geq 3$) 
see the next \S\ref{sec:conj}. 
\end{Rem}

\noindent
{\em Proof of Theorem \ref{thm:main}.} 
(i)$\Rightarrow$(ii): Let $C_0\in\ch^0(\A)$. 
Since 
$$
d_\calL([C_0])=-(q_\infty^{1/2}-q_\infty^{-1/2})[C_0^\lor]+\dots, 
$$
and $q_\infty\neq 1$, we have 
$\rank \left(d_\calL:
\bbC[\ch^0(\A)]\rightarrow\bbC[\ch^1(\A)]\right)=1$ 
(and in particular, 
$H^0(\bbC[\ch^\bullet(\A)], d_\calL)=\Ker (d_\calL: 
\bbC[\ch^0]\rightarrow\bbC[\ch^1])=0$). 

To show that 
\begin{itemize}
\item $H^1(\bbC[\ch^\bullet(\A)], d_\calL)=0$, 
\item $H^2(\bbC[\ch^\bullet(\A)], d_\calL)=\Coker
(d_\calL: 
\bbC[\ch^1]\rightarrow\bbC[\ch^2])$ 
has $\{[C]\}_{C\in\bch(\A)}$ as 
a basis, 
\end{itemize}
it suffices to prove that the induced map 
$$
\overline{d_\calL}:
\bbC[\bch^1(\A)]\longrightarrow
\bbC[\uch^2(\A)]
$$
is surjective (hence bijective). However this 
easily follows from Corollary \ref{cor:nonzero}. 

(ii)$\Rightarrow$(iii) is trivial. 

(iii)$\Rightarrow$(i): Let us assume (iii). 
Since $H^2=\Coker(d_\calL: 
\bbC[\ch^1]\rightarrow\bbC[\ch^2])$, 
the assumption implies that 
the induced map 
\begin{equation}
\label{eq:surj}
\bbC[\ch^1(\A)]\longrightarrow\bbC[\uch^2(\A)]
\mbox{ is surjective}. 
\end{equation}
As in the proof of 
Theorem \ref{thm:det}, $d_\calL$ maps 
$N^1$ to $N^2$. 
Thus the induced map 
$$
\begin{array}{cccc}
\widetilde{d_\calL}: &
W^1\oplus\bbC\cdot[C_0^\lor]&\longrightarrow&
W^2\\
& |\wr & & |\wr\\
&\bbC[\ch^1]/N^1&&\bbC[\uch^2]/N^2
\end{array}
$$
is surjective. Now if $q_\infty=1$, then 
$\widetilde{d_\calL}$ is the zero map 
on $W^1$, and hence $W^2$ is at most one dimensional. 
This is a contradiction to the assumption $\A$ is 
indecomposable (see Proposition \ref{prop:indec}). 
Thus we have $q_\infty\neq 1$. 

Set $\bch^1(\A)=\{C_1, \dots, C_k\}$. Then 
$\ch^1(\A)=\{C_0^\lor, C_1, \dots, C_k\}$. 
$d_\calL([C_0])$ is expressed as 
$$
d_\calL([C_0])=
a_0[C_0^\lor]+
\sum_{i=1}^k a_i[C_i]. 
$$
Note that $a_0=-(q_\infty^{1/2}-q_\infty^{-1/2})\neq 0$. 
Since $d_\calL^2=0$, 
$d_\calL([C_0^\lor])\in\bbC[\ch^2]$ can be expressed 
as a linear combination of 
$d_\calL([C_1]), \dots, d_\calL([C_k])$. The assumption 
(\ref{eq:surj}) implies that (recall that 
$\bbC[\ch^1(\A)]=\bbC[\bch^1(\A)]\oplus\bbC\cdot [C_0^\lor]$)
\begin{equation}
\label{eq:sub}
\bbC[\bch^1(\A)]\longrightarrow\bbC[\uch^2(\A)]
\mbox{ is surjective}. 
\end{equation}
Again by Theorem \ref{thm:det}, we conclude that 
$q_X\neq 1$ for any dense edge $X\in\sfD_\infty(\A_\infty)$ 
in $H_\infty$. 
\qed

\begin{Rem}
\normalfont
\label{rem:counter}
The assumption ``$\A$ is indecomposable'' is 
necessary to prove (iii) $\Rightarrow$ (i) in 
Theorem \ref{thm:main}. 
Indeed, consider the arrangement in 
Figure \ref{fig:buch}, which is decomposable. Let $\calL$ be 
a rank one local system such that $q_1, q_2, q_3\in\bbC^*$ are 
generic and $q_4=q_1^{-1}q_2^{-1}q_3^{-1}$. Then $q_\infty=1$. 
The map 
$\widetilde{d_\calL}:\bbC[\ch^1]\rightarrow \bbC[\uch^2]$ is 
computed as: 
$$
\begin{array}{rcrrr}
\widetilde{d_\calL}([C_1])&=&-(q_{34}^{\frac{1}{2}}-q_{34}^{-\frac{1}{2}})[C_1^\lor]&&\\
\widetilde{d_\calL}([C_2])&=&(q_{234}^{\frac{1}{2}}-q_{234}^{-\frac{1}{2}})[C_1^\lor]&&-(q_{12}^{\frac{1}{2}}-q_{12}^{-\frac{1}{2}})[C_3^\lor]\\
\widetilde{d_\calL}([C_3])&=&&&-(q_{123}^{\frac{1}{2}}-q_{123}^{-\frac{1}{2}})[C_3^\lor]\\
\widetilde{d_\calL}([C_0^\lor])&=&-(q_{2}^{\frac{1}{2}}-q_{2}^{-\frac{1}{2}})[C_1^\lor]&-(q_{12}^{\frac{1}{2}}-q_{12}^{-\frac{1}{2}})[C_2^\lor].&
\end{array}
$$
Hence the map 
$\widetilde{d_\calL}:\bbC[\ch^1]\rightarrow \bbC[\uch^2]$ has 
rank three. This implies that $H^2(\bbC[\ch^\bullet], d_\calL)$ is 
generated by $\bch^2=\{C_4\}$. 
Thus (iii) holds true, while 
(i) is false because $q_\infty=1$. 
\end{Rem}

\section{Remarks and conjectures}
\label{sec:conj}

We conclude this paper with some remarks 
on higher dimensional cases $\ell\geq 3$. 
As in the case $\ell=2$, it seems natural to 
focus on the induced map 
$$
\overline{d_\calL}: \bbC[\bch^{q-1}]\longrightarrow 
\bbC[\uch^{q}]
$$
defined by the composition 
$\bbC[\bch^{q-1}]\hookrightarrow
\bbC[\ch^{q-1}]\stackrel{d_\calL}{\longrightarrow} 
\bbC[\ch^{q}]\twoheadrightarrow\bbC[\uch^{q}]$. 
Since the bases of two spaces 
$\bbC[\bch^{q-1}]$ and $\bbC[\uch^{q}]$ 
are 
naturally identified by the 
involution $\iota$, it makes sense to 
consider the determinant of $\overline{d_\calL}$. 

\begin{Conj}
\normalfont\label{conj:det}
The determinant 
$\det\left(\overline{d_\calL}: \bbC[\bch^{q-1}]\rightarrow 
\bbC[\uch^{q}]\right)$ 
is expressed in the following form 
$$
\det(\overline{d_\calL})=\pm
\prod_{X}(q_X^{1/2}-q_X^{-1/2})^{n_X}, 
$$
where $X$ runs all dense edge $X\in\sfD_\infty(\A_\infty)$ 
with $\dim X\geq \ell-q$ and $n_X>0$. 
\end{Conj}
Once the above conjecture is established, it 
deduces the following. 
\begin{Conj}
\normalfont\label{conj:gen}
Let $\A$ be an essential affine arrangement in $\bbR^\ell$. 
If the rank one local system $\calL$ satisfies the 
condition (\ref{eq:cdo}), then (\ref{eq:gen}) holds. 
\end{Conj}

\noindent
{\em ``Proof of \ref{conj:det} $\Rightarrow$ \ref{conj:gen}.''} 
Since the composition 
$\overline{d_\calL}: \bbC[\bch^{q-1}]\hookrightarrow
\bbC[\ch^{q-1}]\stackrel{d_\calL}{\longrightarrow} 
\bbC[\ch^{q}]\twoheadrightarrow\bbC[\uch^{q}]$ is 
bijective, the rank of the map 
$\bbC[\ch^{q-1}]\stackrel{{d_\calL}}{\longrightarrow} 
\bbC[\ch^{q}]$ is at least $\sharp \bch^{q-1}=
\sharp\uch^{q}$. Hence, 
\begin{eqnarray*}
\dim\Im(d_\calL:\bbC[\ch^{q-1}]\rightarrow\bbC[\ch^{q}])&\geq&
\sharp\uch^{q}=\sharp \bch^{q-1}, \\
\dim\Ker(d_\calL:\bbC[\ch^{q}]\rightarrow\bbC[\ch^{q+1}])&\leq&
\sharp\ch^{q}-
\sharp\bch^{q}=
\sharp\uch^{q}, 
\end{eqnarray*}
for $q\leq \ell-1$. This implies $H^k(\bbC[\ch^\bullet], d_\calL)=0$ 
for $k\leq \ell-1$. Also this deduces 
$H^\ell(\bbC[\ch^\bullet], d_\calL)$ is generated by 
$\bch^\ell(\A)=\bch(\A)$. \qed

\medskip
\noindent
{\bf Acknowledgements}. 
This is an expanded version of the author's talk 
at 
``Fifth Franco-Japanese Symposium on Singularities.'' 
He is grateful to the organizers. 
This work was supported by 
JSPS Grant-in-Aid for Young Scientists (B) 
20740038.

\end{document}